\documentclass[11pt,a4paper]{amsart}
\usepackage[T1]{fontenc}
\usepackage{lmodern,microtype,amsmath,amssymb,amsthm,mathtools,booktabs,array}
\usepackage{tikz}
\usepackage[colorlinks=true,linkcolor=blue!45!black,citecolor=blue!45!black,urlcolor=blue!45!black]{hyperref}
\usepackage[nameinlink,noabbrev]{cleveref}
\newtheorem{theorem}{Theorem}[section]
\newtheorem{lemma}[theorem]{Lemma}
\newtheorem{proposition}[theorem]{Proposition}
\newtheorem{corollary}[theorem]{Corollary}
\theoremstyle{definition}\newtheorem{definition}[theorem]{Definition}
\theoremstyle{remark}
\newcommand{\E}{\mathbb E}
\newcommand{\PP}{\mathbb P}
\newcommand{\Ps}{\mathbb P_\star}
\newcommand{\Es}{\mathbb E_\star}
\newcommand{\Poi}{\operatorname{Poisson}}
\newcommand{\TV}{d_{\mathrm{TV}}}
\newcommand{\ind}{\mathbf 1}
\newcommand{\Ham}{d_{\mathrm H}}
\newcommand{\calL}{\mathcal L}

\DeclareMathOperator{\Var}{Var}
\title[Solution Independence in Random $K$-SAT]{Solution Independence and Poisson Limits in Random $K$-SAT with Growing Width}
\author{Jun Liu}
\address{School of Mathematics, Taiyuan University of Technology, Jinzhong, 030600, China}
\email{liujun03@tyut.edu.cn}
\author{Chunyan Zhao}
\address{School of Data Science and Engineering, East China Normal University, Shanghai, 200062, China}
\email{cyzhao@dase.ecnu.edu.cn}
\date{}
\begin{document}
\begin{abstract}
We study the solution count $X$ of random $K$-SAT when $\E X\to\lambda\in(0,\infty)$. Each clause contains $K$ distinct variables, and clauses are sampled independently. We prove that $X$ converges to $\Poi(\lambda)$ if and only if $N2^{-K}\to0$; the same condition characterizes critical second-factorial-moment independence. At the boundary $N2^{-K}\to a\in(0,\infty)$, the second factorial moment tends to $\lambda^2+\lambda(e^a-1)$. Conditioning on a fixed solution reveals the source of this correction: the number of free coordinate directions converges to $\Poi(a)$, and all their combinations remain solutions with high probability. If $N2^{-K}\to\infty$, satisfiability has probability tending to zero despite the positive limiting mean. For $0<\lambda<1$, unique satisfiability and unsatisfiability coexist with positive lower limiting probabilities if and only if $N2^{-K}$ is eventually bounded. Thus coexistence holds throughout a width range where the ordinary Poisson limit fails. In this range, almost all clauses of a uniquely satisfiable formula are redundant with high conditional probability. This gives support-preserving one-clause solution blocking with $\Theta(N^2)$ clauses at bounded offsets from $\log_2N$.

\end{abstract}
\maketitle

\section{Introduction}
How large must the clause width be for the solutions of random $K$-SAT to behave independently? For fixed width, nearby assignments often satisfy or violate the same clauses. Increasing the width reduces this dependence, but the number of nearby assignments also grows with the number of variables. The resulting contribution to the solution count depends on both effects.

We consider formulas with $N$ variables and $M$ independent random clauses, choosing $M$ so that the solution count $X$ satisfies $\E X\to\lambda\in(0,\infty)$. The controlling parameter is
\[
 a_N=N2^{-K}.
\]
When $K\to\infty$ and $K=O(\log N)$, this is asymptotically the expected number of solutions at Hamming distance one from a fixed assignment conditioned to satisfy the formula. We determine how this local dependence affects Poisson counting and the coexistence of uniquely satisfiable and unsatisfiable formulas.

\subsection{Background and related work}
Frieze and Wormald~\cite{FW} located the satisfiability transition for moderately growing width under $K-\log_2N\to\infty$. With
\[
 M_0=\frac{N\ln2}{-\ln(1-2^{-K})},
\]
the transition is bracketed by $(1\pm\epsilon_N)M_0$ whenever $N\epsilon_N\to\infty$. Coja-Oghlan and Frieze~\cite{CF} studied the finer scale $M=2^K(N\ln2+c)$ and obtained the limiting satisfiability probability $1-\exp(-e^{-c})$. Their proof establishes a Poisson limit for the number of solutions. Their stated theorem assumes $K-\log_2N\to\infty$ and $K-\log_2N=o(\ln N)$, and their conclusion discusses extension to larger widths. The model definition in~\cite{CF} samples literals with replacement within a clause. We work throughout with distinct variables within each clause and retain the corresponding hypergeometric pair probability.

Liu, Gao, and Xu~\cite{LGX} lowered the width condition for threshold localization to $K\ge(1/2+\theta)\log_2N$, for fixed $\theta>0$. Their theorem brackets the transition at $(1\pm c/\sqrt N)M_0$, with fixed $c>0$. At the lower endpoint of that window, the expected number of solutions is $2^{c\sqrt N}$, up to rounding effects. Here we study the finer window in which the expected count has a positive finite limit.

The self-referential solution-blocking construction of Fang et al.~\cite{FLL} uses uniquely satisfiable random formulas with redundant clauses. Replacing one redundant clause by a clause falsified by the unique solution produces an unsatisfiable formula on the same unsigned incidence graph. This motivates a separate probabilistic question: can unique satisfiability and unsatisfiability coexist with positive probability even when the Poisson solution-count limit fails? Our boundary estimates answer this question affirmatively in the native random-clause model.

Our contribution is the width criterion together with the boundary analysis and the weaker criterion for coexistence. For Poisson sufficiency, we use positive-relation Poisson approximation~\cite{BHJ,Ross} to give a direct proof for clauses with distinct variables.

\subsection{Results and proof strategy}
The main results distinguish two width requirements. Critical second-factorial-moment independence and convergence to $\Poi(\lambda)$ each hold exactly when $a_N\to0$, or $K-\log_2N\to\infty$. For $0<\lambda<1$, coexistence holds exactly when $a_N=O(1)$, or $K\ge\log_2N-O(1)$.

At the boundary $a_N\to a\in(0,\infty)$, we obtain
\[
 \E[X(X-1)]\longrightarrow\lambda^2+\lambda(e^a-1).
\]
The correction comes from solutions at bounded Hamming distance. Under conditioning on a fixed solution, the number of free coordinate directions converges to $\Poi(a)$, and their combinations form a cube of solutions with high probability. A size-bias identity converts this local structure into an obstruction to the ordinary Poisson limit. Below this width range, $a_N\to\infty$ implies that satisfiability has probability tending to zero.

Section~\ref{sec:pairs} establishes the required Hamming-shell estimates, including uniform control of the distances between the local and central shells. Section~\ref{sec:geometry} proves the local cube theorem and the Poisson criterion. For coexistence, Section~\ref{sec:coexist} estimates the probability of having no free directions after conditioning on two solutions at a central distance. The estimate is uniform in that distance, allowing it to be summed over all central shells. The resulting uniqueness bound also yields the redundant clauses needed for solution blocking.

\section{Model and main results}\label{sec:model}
\subsection{The ensemble and critical normalization}
Let $2\le K=K(N)\le N$. A clause is obtained by choosing a uniform $K$-element subset of $[N]$ and a uniform sign pattern on that subset. A formula $F$ consists of $M=M(N)$ independent clauses. Repeated clauses are permitted and are kept in indexed slots. Let $X$ count satisfying assignments, and put
\begin{equation}\label{eq:parameters}
 p=2^{-K},\quad q=1-p,\quad L=-\ln q,\quad
 \lambda_N=2^Nq^M,\quad a_N=Np.
\end{equation}
Throughout the main results we assume
\begin{equation}\label{eq:critical}
 \lambda_N\longrightarrow\lambda\in(0,\infty).
\end{equation}
The index $N$ may run through any increasing sequence of integers. This convention handles subsequences caused by integer widths and clause counts. When $K\to\infty$, taking
\[
 M=\left\lfloor\frac{N\ln2-\ln\lambda}{L}\right\rfloor
\]
ensures~\eqref{eq:critical}. For bounded $K$, the results apply to sequences satisfying~\eqref{eq:critical}; such a limit need not be attainable for every integer $N$ with an arbitrarily prescribed $\lambda$.

\begin{definition}
The ensemble is \emph{critically second-moment independent} if
\[
 R_N:=\frac{\E[X(X-1)]}{\lambda_N^2}\longrightarrow1.
\]
Unique satisfiability and unsatisfiability \emph{coexist} if
\[
 \liminf_{N\to\infty}\PP(X=1)>0,
 \qquad \liminf_{N\to\infty}\PP(X=0)>0.
\]
\end{definition}
Second-moment independence averages joint survival probabilities over distinct candidate assignments; it imposes no uniform bound on individual pair correlations. We prove this property and Poisson convergence separately.

Fix $\sigma\in\{0,1\}^N$. Write $\Ps$ for the formula law conditional on $\sigma\models F$, and $\Es$ for its expectation. Let $U=U(F,\sigma)$ be the set of coordinates $i$ for which flipping only coordinate $i$ preserves satisfiability, and let $D=|U|$. Write $\sigma^S$ for $\sigma$ with the coordinates of $S$ flipped. This conditioned law is also the usual planted law: each clause is independently uniform among clauses satisfied by $\sigma$.

\subsection{Statements}
\begin{theorem}[Width criteria]\label{thm:criterion}
Under~\eqref{eq:critical}, the following are equivalent:
\begin{enumerate}
\item $a_N\to0$, equivalently $K-\log_2N\to+\infty$;
\item $R_N\to1$;
\item $X\Rightarrow\Poi(\lambda)$.
\end{enumerate}
In the third assertion convergence also holds in total variation.
\end{theorem}

\begin{theorem}[Finite boundary]\label{thm:boundary}
Suppose~\eqref{eq:critical} holds and $a_N\to a\in(0,\infty)$, equivalently $K-\log_2N\to-\log_2a$. Then
\begin{align}
 \E[X(X-1)]&\longrightarrow\lambda^2+\lambda(e^a-1),\label{eq:boundarymoment}\\
 \Var(X)&\longrightarrow\lambda e^a.\label{eq:boundaryvar}
\end{align}
Under $\Ps$, $D\Rightarrow\Poi(a)$ and
\begin{equation}\label{eq:cube}
 \Ps(\sigma^S\models F\text{ for every }S\subseteq U)\longrightarrow1.
\end{equation}
The unconditional solution count does not converge to $\Poi(\lambda)$.
\end{theorem}

\begin{theorem}[Diverging neighbor scale]\label{thm:collapse}
If~\eqref{eq:critical} holds and $a_N\to\infty$, equivalently $K-\log_2N\to-\infty$, then
\[
 R_N\to\infty,\qquad \PP(X>0)\to0.
\]
Thus $X\to0$ in probability although $\E X\to\lambda>0$.
\end{theorem}

\begin{theorem}[Coexistence]\label{thm:coexist}
Assume~\eqref{eq:critical}. If $a_N\to a\in(0,\infty)$, equivalently $K-\log_2N\to-\log_2a$, then
\begin{align}
 \liminf\PP(X=1)&\ge\lambda e^{-a}(1-\lambda),\label{eq:unique-lower}\\
 \limsup\PP(X=1)&\le\lambda e^{-a},\label{eq:unique-upper}\\
 \liminf\PP(X=0)&\ge1-\lambda.\label{eq:zero-lower}
\end{align}
For fixed $0<\lambda<1$, coexistence holds if and only if $\limsup a_N<\infty$, equivalently $K\ge\log_2N-O(1)$. If $A=\limsup a_N<\infty$, then
\[
 \liminf\PP(X=1)\ge\lambda e^{-A}(1-\lambda).
\]
\end{theorem}
The lower bounds give coexistence for $0<\lambda<1$. The full boundary distribution and the sizes of maximal solution components are discussed in Section~\ref{sec:discussion}.

\begin{table}[t]
\centering\small
\setlength{\tabcolsep}{4pt}
\begin{tabular}{@{}>{\raggedright\arraybackslash}p{.29\textwidth}>{\raggedright\arraybackslash}p{.27\textwidth}>{\raggedright\arraybackslash}p{.37\textwidth}@{}}
\toprule
Clause width $K$ & Second factorial moment & Counting and coexistence \\\midrule
$K-\log_2N\to+\infty$ & $\lambda^2$ & $\Poi(\lambda)$; coexistence for all $\lambda>0$ \\[3pt]
$K-\log_2N\to w\in\mathbb R$ & $\lambda^2+\lambda(e^{2^{-w}}-1)$ & Non-Poisson; coexistence proved for $0<\lambda<1$ \\[3pt]
$K-\log_2N\to-\infty$ & Diverges & SAT probability tends to zero \\\bottomrule
\end{tabular}
\caption{The width transition under $\E X\to\lambda\in(0,\infty)$, stated in terms of $K$. In the middle row $a=2^{-w}$; the statement applies along convergent subsequences when integer rounding makes the width offset oscillate.}
\label{tab:regimes}
\end{table}

\section{Pair correlations across Hamming distances}\label{sec:pairs}
\subsection{Exact identities}
For assignments at distance $h$, define
\begin{equation}\label{eq:rho}
 \rho(h)=\frac{\binom{N-h}{K}}{\binom NK},
\end{equation}
with value zero if $N-h<K$. A clause is falsified by both assignments with probability $p\rho(h)$. Hence, with
\[
 s_h=1-2p+p\rho(h),\qquad
 u_h=\frac{p(1-\rho(h))}{q},
\]
we have
\begin{align}
 \PP(\sigma,\tau\models F)&=s_h^M,\label{eq:joint}\\
 T_N(h):=\Es Z_h&=\binom Nh(1-u_h)^M,\label{eq:shell}\\
 R_N&=2^{-N}\sum_{h=1}^N\binom Nh
 \left(1+\frac{p(\rho(h)-p)}{q^2}\right)^M.\label{eq:ratio}
\end{align}
Here $Z_h$ counts solutions at distance $h$ from $\sigma$, so
\begin{equation}\label{eq:rootedmoment}
 \E[X(X-1)]=\lambda_N\sum_{h=1}^NT_N(h).
\end{equation}
Also,
\begin{equation}\label{eq:rhobound}
 \rho(h)\le(1-h/N)^K\le e^{-Kh/N}.
\end{equation}
At the finite boundary, the subtraction of $p$ in~\eqref{eq:ratio} contributes a nonvanishing term. If $a_N\to a\in(0,\infty)$, expansion of $p/L$ gives
\begin{equation}\label{eq:mpboundary}
 Mp=N\ln2-\ln\lambda_N-\tfrac12a\ln2+o(1),
 \qquad Mp^2\to a\ln2.
\end{equation}
Thus the error $O(Mp^2)$ cannot be discarded there.

\subsection{The fixed-width obstruction}
For fixed width, the failure of independence is stronger than the divergence
deduced from single-coordinate neighbors. The following shell argument also
appears in the fixed-width obstruction of~\cite{FLL}.
\begin{proposition}\label{prop:fixed}
For every fixed $K\ge2$ (hence $a_N=N2^{-K}\to\infty$) under~\eqref{eq:critical}, there is a constant
$c_K'>0$ such that $R_N\ge\exp(c_K'N)$ for all sufficiently large $N$.
\end{proposition}
\begin{proof}
Put $c_K=\ln2/L$. Then $M=c_KN+O(1)$. For $h=\lfloor\beta N\rfloor$,
where $\beta\in(0,1)$ is fixed, the contribution of this shell to
\eqref{eq:ratio} has logarithm $N\Phi_K(\beta)+o(N)$, with
\[
 \Phi_K(\beta)=H(\beta)-\ln2+
 c_K\ln\frac{1-2p+p(1-\beta)^K}{q^2}.
\]
The replacement of the exact hypergeometric factor by $(1-\beta)^K$
has error $O(1/N)$ and changes this logarithm by only $O(1)$.
At $\beta=0$ the expression vanishes, while as $\beta\downarrow0$,
\[
 \Phi_K(\beta)=\beta\ln(1/\beta)
 +\left(1-\frac{c_KpK}{q}\right)\beta+O_K(\beta^2).
\]
Here $H$ denotes binary entropy with natural logarithms.
For a sufficiently small fixed positive $\beta$, the expression is positive.
Keeping this one shell proves the assertion.
\end{proof}

\subsection{A uniform shell decomposition}
\begin{lemma}[Shell limits and uniform tails]\label{lem:shells}
Suppose~\eqref{eq:critical}, $K=O(\log N)$, and $a_N\to a\in[0,\infty)$. Equivalently, $K-\log_2N\to w\in\mathbb R\cup\{+\infty\}$ with $a=2^{-w}$ and $2^{-\infty}=0$. There is a constant $\alpha>0$ such that, for $d_N=N^{-\alpha}$ and
\[
 \mathcal C_N=\{h:|h/N-1/2|\le d_N\},
\]
the following hold:
\begin{align}
 T_N(h)&\to a^h/h! &&\text{for every fixed }h\ge1,\label{eq:fixedshell}\\
 \sum_{1\le h\le N/K^3}T_N(h)&\to e^a-1,\label{eq:localsum}\\
 \sum_{h\in\mathcal C_N}T_N(h)&\to\lambda,\label{eq:centralsum}\\
 \sum_{\substack{h>N/K^3\\h\notin\mathcal C_N}}T_N(h)&\to0.\label{eq:outsidesum}
\end{align}
Moreover, for some constant $B$, uniformly for $1\le h\le N/K^3$,
\begin{equation}\label{eq:factorialdomination}
 T_N(h)\le B^h/h!.
\end{equation}
\end{lemma}
\begin{proof}
Boundedness of $a_N$ implies $K\ge\log_2N-O(1)$ and $p=O(1/N)$. Thus
\[
 Mp=N\ln2+O(1),\quad Mp/q=N\ln2+O(1),\quad Mp/q^2=N\ln2+O(1).
\]
For fixed $h$, expansion of~\eqref{eq:rho} yields
\[
 1-\rho(h)=Kh/N+O_h(K^2/N^2).
\]
Consequently $M\ln(1-u_h)=-Kh\ln2+o(1)$, proving~\eqref{eq:fixedshell}.

For $h\le N/K^3$, use $1-e^{-x}\ge x-x^2/2$ in~\eqref{eq:rhobound}. Since $\ln(1-u_h)\le-u_h$,
\[
 M\ln(1-u_h)
 \le-Kh\ln2+C\left(\frac{Kh}{N}+\frac{K^2h^2}{N}\right)
 \le-Kh\ln2+Ch/K.
\]
Together with $\binom Nh\le N^h/h!$, this proves
\[
 T_N(h)\le\frac{(a_Ne^{C/K})^h}{h!}\le\frac{B^h}{h!}.
\]
The summable majorant proves~\eqref{eq:localsum} and also the uniform tail statement
\begin{equation}\label{eq:uniformtail}
 \lim_{H\to\infty}\limsup_{N\to\infty}
 \sum_{H<h\le N/K^3}T_N(h)=0.
\end{equation}

For $N/K^3<h\le N/K$, the inequality $1-e^{-x}\ge(1-e^{-1})x$ for $0\le x\le1$ gives
\[
 T_N(h)\le\exp\{h\ln(eN/h)-cKh\}.
\]
Here $\ln(eN/h)\le1+3\ln K=o(K)$, so the sum is at most $N\exp(-c'N/K^2)=o(1)$.

Choose $b\in(0,1/2)$ so small that $H(b)<(\ln2)(1-e^{-1})$, where $H(x)=-x\ln x-(1-x)\ln(1-x)$. For $N/K<h\le bN$, we have $\rho(h)\le e^{-1}$. The entropy bound for binomial coefficients shows that the sum in this interval is at most
\[
 (N+1)\exp\{NH(b)-(N\ln2+O(1))(1-e^{-1})\}=e^{-\Omega(N)}.
\]

For the remaining lower tail, put $\beta=-\log_2(1-b)>0$ and choose $0<\alpha<\min\{1/2,\beta/2\}$. If $h>bN$, then $\rho(h)\le C N^{-\beta}$. From~\eqref{eq:ratio} and~\eqref{eq:rootedmoment},
\begin{equation}\label{eq:shellratio}
 T_N(h)=\lambda_N\,2^{-N}\binom Nh
 \left(1+\frac{p(\rho(h)-p)}{q^2}\right)^M.
\end{equation}
The last factor is at most $\exp(O(N^{1-\beta}))$. A binomial tail bound therefore gives
\[
 \sum_{bN<h<(1/2-d_N)N}T_N(h)
 \le O(1)\exp\{O(N^{1-\beta})-2N^{1-2\alpha}\}=o(1).
\]

Uniformly on $\mathcal C_N$, the product in~\eqref{eq:rho} gives
\begin{equation}\label{eq:centralrho}
 \rho(h)/p=1+O(Kd_N+K^2/N)=1+o(1).
\end{equation}
Since $Mp^2=O(1)$, the logarithm of the last factor in~\eqref{eq:shellratio} is $o(1)$ uniformly. The binomial mass of $\mathcal C_N$ tends to one, proving~\eqref{eq:centralsum}. Finally, if $h>(1/2+d_N)N$, then $\rho(h)\le p$, so that factor is at most one. The upper tail contributes $o(1)$. This proves~\eqref{eq:outsidesum}.
\end{proof}

\begin{corollary}\label{cor:boundedmoments}
Under the assumptions of \cref{lem:shells},
\[
 \E[X(X-1)]\to\lambda^2+\lambda(e^a-1),
 \qquad \Var(X)\to\lambda e^a.
\]
\end{corollary}
\begin{proof}
Sum the shell limits and use~\eqref{eq:rootedmoment}. Then use $\E X=\lambda_N\to\lambda$.
\end{proof}

\subsection{A collision bound above the boundary}
\begin{lemma}\label{lem:collision}
If $a_N\to0$, equivalently $K-\log_2N\to+\infty$, and $C_N=N\ln2+O(1)$, then
\begin{equation}\label{eq:collision}
 2^{-N}\sum_{h=1}^N\binom Nh e^{C_N\rho(h)}=1+o(1).
\end{equation}
This statement holds for all $K\le N$.
\end{lemma}
\begin{proof}
First suppose $K=O(\log N)$. For $h\le N/K^3$, the proof of~\eqref{eq:factorialdomination} now gives
\[
 2^{-N}\binom Nh e^{C_N\rho(h)}
 \le C\frac{(a_Ne^{C'/K})^h}{h!},
\]
whose sum tends to zero. The next three intervals in the proof of \cref{lem:shells} have the same vanishing bounds. On $h\ge(1/2-d_N)N$, we have $\rho(h)\le p(1+o(1))$, so $C_N\rho(h)=o(1)$ uniformly. That interval has binomial mass tending to one. This proves~\eqref{eq:collision}.

For general $K$, let $k=\min\{K,\lceil2\log_2N\rceil\}$. Then $k=O(\log N)$ and $N2^{-k}\to0$. The probability $\rho_{N,K}(h)$ is nonincreasing in $K$, so the sum for $K$ is bounded above by the sum for $k$, once $C_N>0$. Its lower bound is $1-2^{-N}$. The result follows.
\end{proof}

\section{Local geometry and the Poisson criterion}\label{sec:geometry}
\subsection{Size bias and free directions}
For any bounded function $f$ of the solution count, symmetry gives
\begin{equation}\label{eq:sizebias}
 \Es f(X)=\frac{\E[Xf(X)]}{\lambda_N}.
\end{equation}
Indeed, sum $\E[\ind_{\{\sigma\models F\}}f(X)]$ over all $\sigma$. This identity concerns the conditional law of the count, not a uniform choice of a satisfiable formula.

A clause satisfied by $\sigma$ blocks direction $i$ exactly when its only true literal under $\sigma$ lies on variable $i$. Under $\Ps$ the probability of this category is
\begin{equation}\label{eq:blocker}
 b=\frac{pK}{Nq}.
\end{equation}
The categories for different $i$ are disjoint. If $D$ counts categories absent from all $M$ clauses, then
\begin{equation}\label{eq:freefactorial}
 \Es(D)_r=(N)_r(1-rb)^M\qquad(0\le r\le N).
\end{equation}

\begin{lemma}\label{lem:free}
If $a_N\to a\in(0,\infty)$, equivalently $K-\log_2N\to-\log_2a$, then under $\Ps$,
\[
 D\Rightarrow\Poi(a),\qquad \Ps(D=0)\to e^{-a}.
\]
There is a constant $B$ such that $\Es(D)_r\le B^r$ for all $r\ge0$ and all sufficiently large $N$.
\end{lemma}
\begin{proof}
In this regime $Mb=K\ln2+o(1)$ and $Mb^2=o(1)$. For fixed $r$,~\eqref{eq:freefactorial} tends to $a^r$. Moreover,
\[
 \Es(D)_r\le (Ne^{-Mb})^r\le B^r,
\]
with $(D)_r=0$ for $r>N$. In the expansion
\[
 \Es z^D=\sum_{r\ge0}\frac{(z-1)^r}{r!}\Es(D)_r,
\]
the bound permits dominated convergence for every fixed $z\in[0,1]$. The limit is $e^{a(z-1)}$, proving the assertions.
\end{proof}

\begin{lemma}[Local compatibility]\label{lem:compatibility}
Suppose $a_N\to a\in(0,\infty)$, equivalently $K-\log_2N\to-\log_2a$. For a fixed set $S$ let $A_S=\{S\subseteq U\}$ and $B_S=\{\sigma^S\models F\}$. For each fixed $h\ge1$,
\begin{equation}\label{eq:mismatch}
 \sum_{|S|=h}\Ps(A_S\mathbin\triangle B_S)=O_h(K^2/N)=o(1).
\end{equation}
Consequently~\eqref{eq:cube} holds.
\end{lemma}
\begin{proof}
Work with one clause under the planted law. Let $\mathcal U_S$ be the event that the clause blocks one of the directions in $S$, and let $\mathcal V_S$ be the event that it is falsified by $\sigma^S$. Their probabilities are $hb$ and $u_h$, respectively.

If $\mathcal V_S\setminus\mathcal U_S$ occurs, all true literals under $\sigma$ lie in $S$, and there are at least two. For a fixed nonempty $T\subseteq S$, the probability that precisely the variables in $T$ carry true literals is
\[
 \frac pq\frac{(K)_{|T|}}{(N)_{|T|}}.
\]
Summing over $|T|\ge2$ bounds this event by $O_h(pK^2/N^2)$. If $\mathcal U_S\setminus\mathcal V_S$ occurs, a clause with sole true literal on some $i\in S$ must also use another coordinate of $S$. Its probability is at most
\[
 hb\,\frac{(K-1)(h-1)}{N-1}=O_h(pK^2/N^2).
\]

For $A_S\setminus B_S$, choose a clause witnessing $\mathcal V_S\setminus\mathcal U_S$; every other clause avoids $\mathcal U_S$. For $B_S\setminus A_S$, reverse the roles. Independence of the clauses gives
\begin{align*}
 \Ps(A_S\setminus B_S)&\le M\,O_h(pK^2/N^2)(1-hb)^{M-1},\\
 \Ps(B_S\setminus A_S)&\le M\,O_h(pK^2/N^2)(1-u_h)^{M-1}.
\end{align*}
For fixed $h$, both $\binom Nh(1-hb)^{M-1}$ and $\binom Nh(1-u_h)^{M-1}$ are bounded. Since $Mp=O(N)$, summing proves~\eqref{eq:mismatch}.

Fix an integer $H$. If $D\le H$ and no mismatch occurs for any $S$ of size at most $H$, every subset of $U$ can be flipped. Hence the failure probability in~\eqref{eq:cube} is at most
\[
 \Ps(D>H)+\sum_{h=1}^H\sum_{|S|=h}\Ps(A_S\mathbin\triangle B_S).
\]
Let $N\to\infty$ and then $H\to\infty$, using \cref{lem:free}.
\end{proof}

\begin{figure}[t]
\centering
\begin{tikzpicture}[scale=.9,every node/.style={font=\small}]
\fill (0,0) circle (2.5pt);
\node at (0,-.65) {$D=0$: one solution};
\draw[thick] (4,0)--(5.2,0);
\fill (4,0) circle (2.5pt);\fill (5.2,0) circle (2.5pt);
\node at (4.6,-.65) {$D=1$: two solutions};
\draw[thick] (9,0)--(10.2,0)--(10.2,1.2)--(9,1.2)--cycle;
\foreach \x in {9,10.2}\foreach \y in {0,1.2}\fill (\x,\y) circle (2.5pt);
\node at (9.6,-.65) {$D=2$: four solutions};
\end{tikzpicture}
\caption{Local solutions generated by free coordinate directions. Under conditioning on a fixed solution and $a_N\to a>0$, the number of directions tends to $\Poi(a)$. The diagram depicts guaranteed cubes; it does not identify maximal components or their unconditional frequencies.}
\label{fig:cubes}
\end{figure}

\subsection{The boundary obstruction and low-width collapse}
\begin{proposition}\label{prop:nonpoisson}
If $a_N\to a\in(0,\infty)$, equivalently $K-\log_2N\to-\log_2a$, then $X$ does not converge to $\Poi(\lambda)$.
\end{proposition}
\begin{proof}
Suppose otherwise. For each $j\ge1$,~\eqref{eq:sizebias} gives
\[
 \Ps(X=j)=\frac{j\PP(X=j)}{\lambda_N}
 \longrightarrow e^{-\lambda}\frac{\lambda^{j-1}}{(j-1)!}.
\]
The limiting masses sum to one, so under $\Ps$ we would have $X\Rightarrow1+\Poi(\lambda)$. On the event in~\eqref{eq:cube}, $X\ge2^D$. Thus for each fixed integer $d\ge1$,
\[
 \PP(1+\Poi(\lambda)\ge2^d)
 \ge\PP(\Poi(a)\ge d).
\]
For large fixed $d$ this is impossible: the right side is at least $e^{-a}a^d/d!$, while the left side is at most $(e\lambda/(2^d-1))^{2^d-1}$ once $2^d-1>\lambda$. Their logarithms have orders $-d\ln d$ and $-d2^d$, respectively. This contradiction proves the proposition.
\end{proof}

\begin{proof}[Proof of \cref{thm:collapse}]
Let $\mu_D=\Es D=N(1-b)^M$. Since the blocking categories are disjoint,
\[
 \Es D(D-1)=N(N-1)(1-2b)^M\le N(N-1)(1-b)^{2M},
\]
so $\Var_\star(D)\le\mu_D$. We claim $\mu_D\to\infty$ when $a_N\to\infty$.

Along a subsequence with $K\to\infty$, we have $K\le\log_2N$ eventually, and
\[
 Mb=\frac{K\ln2\,p}{qL}+O(K/N)
 =K\ln2+O(Kp+K/N),\qquad Mb^2=o(1).
\]
Because $Kp\to0$, this gives $\mu_D\sim N2^{-K}=a_N\to\infty$. Along a subsequence with fixed $K$, $Mb$ tends to a finite constant, so $\mu_D=\Theta(N)$. A subsequence argument establishes the claim for arbitrary $K$ with $a_N\to\infty$.

Chebyshev's inequality now implies $D\to\infty$ in planted probability. Since $X\ge D+1$,
\[
 \PP(X>0)=\lambda_N\Es\frac1X
 \le\lambda_N\Es\frac1{D+1}\to0.
\]
Also $\Es(X-1)\ge\mu_D\to\infty$, so~\eqref{eq:rootedmoment} gives $R_N\to\infty$.
\end{proof}

\subsection{Poisson sufficiency and completion of the criterion}
\begin{proposition}\label{prop:poisson}
If $a_N\to0$, equivalently $K-\log_2N\to+\infty$, then $R_N\to1$ and $\TV(\calL(X),\Poi(\lambda))\to0$.
\end{proposition}
\begin{proof}
Here $Mp=N\ln2-\ln\lambda_N+o(1)$ and $Mp^2=o(1)$. Uniformly in $h$, the logarithm of the factor in~\eqref{eq:ratio} is $Mp\rho(h)+O(Mp^2)$. \Cref{lem:collision} proves $R_N\to1$.

For the distributional assertion, first take a Poisson number of clauses with mean $m=M+o(1/p)$, and write $Y_m$ for the solution count. The multiplicities of individual clauses are independent Poisson variables. A fixed assignment survives with probability $\pi=e^{-mp}$, and $\mu=\E Y_m=2^Ne^{-mp}\to\lambda$. Conditioning that assignment to survive is achieved by setting to zero the multiplicities of all clauses that falsify it. This coupling only increases the other survival indicators, so they are positively related.

For a sum $Y$ of positively related indicators with means $\pi_i$ and total mean $\mu$, the standard Stein--Chen bound~\cite{BHJ,Ross} gives
\begin{equation}\label{eq:stein}
 \TV(\calL(Y),\Poi(\mu))
 \le\min\{1,\mu^{-1}\}
 \left(\Var(Y)-\mu+2\sum_i\pi_i^2\right).
\end{equation}
For assignments at distance $h$, the joint survival probability is $e^{-mp(2-\rho(h))}$. Hence \cref{lem:collision} implies
\[
 \frac{\E[Y_m(Y_m-1)]}{\mu^2}
 =2^{-N}\sum_{h=1}^N\binom Nh e^{mp\rho(h)}\to1.
\]
Also $\sum_i\pi_i^2=2^{-N}\mu^2\to0$. The right side of~\eqref{eq:stein} tends to zero.

Choose $b_N\to\infty$ such that $b_Np\sqrt M\to0$ and $b_N/\sqrt M\to0$, which is possible because $Mp^2\to0$. Let $m_\pm=M\pm b_N\sqrt M$. Using a common infinite sequence of random clauses and independent Poisson counts $L_\pm$ of these means, Chebyshev's inequality gives
\[
 \PP(L_-\le M\le L_+)\to1.
\]
On this event $Y_{m_+}\le X\le Y_{m_-}$. Both bounding counts converge to $\Poi(\lambda)$ because $p|m_\pm-M|=o(1)$. Their distribution functions sandwich that of $X$ at every half-integer. Thus $X\Rightarrow\Poi(\lambda)$. For integer-valued laws, convergence of point masses to a probability law implies total variation convergence by truncation to a finite set.
\end{proof}

\begin{proof}[Completion of \cref{thm:criterion,thm:boundary}]
\Cref{cor:boundedmoments,lem:free,lem:compatibility,prop:nonpoisson} prove \cref{thm:boundary}. Sufficiency in \cref{thm:criterion} is \cref{prop:poisson}. If $a_N$ does not tend to zero, a subsequence has $a_N\to a\in(0,\infty]$. For finite $a$, \cref{thm:boundary} rules out both independence and the Poisson limit. For $a=\infty$, \cref{thm:collapse} does so. This proves necessity.
\end{proof}

\section{Coexistence beyond solution independence}\label{sec:coexist}
\subsection{Free directions conditional on two distant solutions}
\begin{lemma}[Uniform double conditioning]\label{lem:double}
Suppose $a_N\to a\in(0,\infty)$, equivalently $K-\log_2N\to-\log_2a$, and let $\mathcal C_N$ be the central distance interval of \cref{lem:shells}. Then
\begin{equation}\label{eq:doublezero}
 \sup_{\Ham(\sigma,\tau)\in\mathcal C_N}
 \left|\PP(D=0\mid\sigma,\tau\models F)-e^{-a}\right|\to0,
\end{equation}
where $D$ counts free directions of $\sigma$.
\end{lemma}
\begin{proof}
Conditioning two assignments to satisfy the formula leaves clauses independent, each uniform among clauses satisfied by both. For $h=\Ham(\sigma,\tau)$, put
\[
 \eta_h=\frac{\binom{N-h}{K-1}}{\binom{N-1}{K-1}}.
\]
A clause blocks the direction $i$ at $\sigma$ with conditional probability
\begin{equation}\label{eq:doubleblock}
 b_i^{(h)}=\frac{pK}{Ns_h}
 \begin{cases}
 1,&\sigma_i=\tau_i,\\
 1-\eta_h,&\sigma_i\ne\tau_i.
 \end{cases}
\end{equation}
For an agreeing coordinate, the sole true literal at $i$ satisfies the clause under both assignments. For a disagreeing coordinate, the other $K-1$ variables must not all lie in their agreement set. This gives~\eqref{eq:doubleblock} exactly.

Uniformly for $h\in\mathcal C_N$, $\eta_h=O(p)$ and $s_h=1+O(p)$, whence
\begin{equation}\label{eq:uniformblock}
 Mb_i^{(h)}=K\ln2+o(1)
\end{equation}
uniformly in both $h$ and $i$. Also $M\max_i(b_i^{(h)})^2=O(pK^2/N)=o(1)$. The blocking categories remain disjoint, so for every fixed $r$,
\[
 \E[(D)_r\mid\sigma,\tau\models F]
 =\sum_{i_1,\ldots,i_r\ \mathrm{distinct}}
 \left(1-\sum_{j=1}^r b_{i_j}^{(h)}\right)^M
 \longrightarrow a^r
\]
uniformly over central pairs.

Let $b_{\min}$ be the minimum of $b_i^{(h)}$ over these pairs and coordinates. By~\eqref{eq:uniformblock}, $Ne^{-Mb_{\min}}$ is bounded. For all $r$, the displayed factorial moment is at most $(Ne^{-Mb_{\min}})^r\le B^r$. Inclusion--exclusion therefore gives an absolutely dominated expansion
\[
 \PP(D=0\mid\sigma,\tau\models F)
 =\sum_{r\ge0}\frac{(-1)^r}{r!}
 \E[(D)_r\mid\sigma,\tau\models F].
\]
Uniform convergence of each fixed factorial moment and the common summable bound prove~\eqref{eq:doublezero}.
\end{proof}

\subsection{Isolated roots and the coexistence criterion}
\begin{lemma}\label{lem:isolated}
If $a_N\to a\in(0,\infty)$, equivalently $K-\log_2N\to-\log_2a$, then
\begin{equation}\label{eq:isolatedmean}
 \Es[(X-1)\ind_{\{D=0\}}]\longrightarrow\lambda e^{-a}.
\end{equation}
\end{lemma}
\begin{proof}
For fixed $h\ge1$, $D=0$ implies $A_S$ fails for every nonempty $S$. Thus
\[
 \Es[Z_h\ind_{\{D=0\}}]
 \le\sum_{|S|=h}\Ps(B_S\setminus A_S)=o(1)
\]
by \cref{lem:compatibility}. For $H<h\le N/K^3$, dominate by $T_N(h)$ and use~\eqref{eq:uniformtail}. All noncentral distances larger than $N/K^3$ contribute $o(1)$ by~\eqref{eq:outsidesum}.

For central distances, condition on the second assignment being a solution. \Cref{lem:double} gives the uniform multiplier $e^{-a}+o(1)$. Therefore their sum is
\[
 (e^{-a}+o(1))\sum_{h\in\mathcal C_N}T_N(h)
 \longrightarrow\lambda e^{-a}.
\]
Letting first $N\to\infty$ and then $H\to\infty$ completes the proof.
\end{proof}

\begin{proof}[Proof of \cref{thm:coexist}]
Under $\Ps$, $X\ge1$. If $X=1$ then $D=0$, and the pointwise inequalities
\[
 \ind_{\{D=0\}}-(X-1)\ind_{\{D=0\}}
 \le\ind_{\{X=1\}}\le\ind_{\{D=0\}}
\]
hold. Apply \cref{lem:free,lem:isolated}, and use $\PP(X=1)=\lambda_N\Ps(X=1)$ from~\eqref{eq:sizebias}. This proves~\eqref{eq:unique-lower} and~\eqref{eq:unique-upper}. Markov's inequality yields~\eqref{eq:zero-lower}.

Suppose now $0<\lambda<1$ and $A=\limsup a_N<\infty$. Any subsequence has a further subsequence with $a_N\to a\in[0,A]$. If $a>0$, the lower bound just proved is at least $\lambda e^{-A}(1-\lambda)$. If $a=0$, \cref{prop:poisson} gives $\PP(X=1)\to\lambda e^{-\lambda}\ge\lambda(1-\lambda)$, which is at least the same bound. Hence the positive lower bound holds along the full sequence, and unsatisfiability has lower limit at least $1-\lambda>0$.

Conversely, if $\limsup a_N=\infty$, some subsequence has $a_N\to\infty$. By \cref{thm:collapse}, $\PP(X=1)\le\PP(X>0)\to0$ along it, ruling out coexistence.
\end{proof}

\begin{corollary}[Quadratic clause counts]\label{cor:quadratic}
Let $K=\log_2N+O(1)$ be integer-valued, equivalently $0<\liminf a_N\le\limsup a_N<\infty$, and choose $M$ so that $\lambda_N\to\lambda\in(0,1)$. Then $M=\Theta(N^2)$, unique satisfiability and unsatisfiability coexist, and $X$ does not converge to $\Poi(\lambda)$.
\end{corollary}
\begin{proof}
The sequence $a_N$ is bounded above and bounded away from zero. Apply \cref{thm:criterion,thm:coexist}, and use $M\sim2^KN\ln2$.
\end{proof}

\subsection{Redundant clauses and support-preserving solution blocking}
For the present critical ensemble, the redundancy conclusion of \cite[Lemma~5]{FLL} holds whenever unique satisfiability has probability bounded away from zero. The proof counts assignments that violate exactly one clause; Poisson convergence is unnecessary.

\begin{theorem}[Redundancy under positive-probability uniqueness]\label{thm:redundancy}
Assume~\eqref{eq:critical} and $\liminf\PP(X=1)>0$. Then $a_N=O(1)$, equivalently $K\ge\log_2N-O(1)$, and $M=\Omega(N^2)$. Let $R(F)$ be the number of redundant clause slots. Then
\begin{align}
 \PP\bigl(R(F)\ge M-N\ln N\mid X=1\bigr)&=1-o(1),\label{eq:conditionalredundancy}\\
 \PP\bigl(X=1,\ R(F)\ge M-N\ln N\bigr)&=\PP(X=1)-o(1).\label{eq:jointredundancy}
\end{align}
In particular, conditional on unique satisfiability, a $1-o(1)$ fraction of clause slots are redundant with probability $1-o(1)$.
\end{theorem}
\begin{proof}
If $a_N$ were unbounded, a subsequence would have $a_N\to\infty$. By \cref{thm:collapse}, unique satisfiability would have probability tending to zero along it, contradicting the assumed positive lower limit. Thus $a_N=O(1)$, and $M\sim 2^KN\ln2=\Omega(N^2)$.

Let $W$ count assignments violating exactly one indexed clause. Then
\[
 \E W=2^NMpq^{M-1}=\lambda_N\frac{Mp}{q}=\Theta(N).
\]
Markov's inequality gives $\PP(W>N\ln N)=O(1/\ln N)=o(1)$. If $F$ is uniquely satisfiable, every essential slot has an assignment that satisfies all other slots but violates that slot. Different essential slots require different such assignments. Hence the number of essential slots is at most $W$. On $\{X=1,W\le N\ln N\}$, at least $M-N\ln N$ slots are redundant. This proves~\eqref{eq:jointredundancy}; dividing the exceptional probability by $\PP(X=1)$ proves~\eqref{eq:conditionalredundancy}. Finally $N\ln N=o(M)$.
\end{proof}

Coexistence therefore suffices for redundancy in this model, although only the positive lower bound on $\PP(X=1)$ enters the proof. Equation~\eqref{eq:jointredundancy} also identifies the relevant probability: if $\PP(X=1)\to\pi_1>0$, then the joint probability tends to $\pi_1$. This recovers $\lambda e^{-\lambda}$ in the Poisson regime of~\cite{FLL}; at the finite boundary, \cref{thm:coexist} supplies a positive lower limit for $0<\lambda<1$. The argument uses both the random-clause model and normalization~\eqref{eq:critical}.

\begin{corollary}\label{cor:blocking}
Under the assumptions of \cref{cor:quadratic}, with probability bounded away from zero the formula is uniquely satisfiable and has at least $M-N\ln N$ redundant clause slots. Each such formula admits replacement of one clause by a $K$-clause that makes the formula unsatisfiable and preserves its unsigned variable--clause incidence graph.
\end{corollary}
\begin{proof}
By \cref{cor:quadratic,thm:redundancy}, the event of unique satisfiability with at least $M-N\ln N$ redundant slots has probability bounded away from zero.

Choose a redundant slot with variable set $J$ and let $\sigma$ be the unique solution. Replace its clause by the clause on $J$ whose every literal is false under $\sigma$. Deleting the original redundant clause leaves exactly the solution $\sigma$; inserting the new clause removes it. The replacement uses the same variables and preserves every unsigned incidence. Since $N\ln N=o(M)$, there are redundant slots for all sufficiently large $N$ on the stated event.
\end{proof}
This extends the probabilistic basis of the solution-blocking construction in~\cite{FLL} to bounded width offsets, where Poisson convergence fails. The construction is existential: finding the unique solution and certifying redundancy remain computational tasks.

\section{Discussion}\label{sec:discussion}
\subsection{Width and the critical window}
The distinction between threshold localization and critical counting is essential when comparing width conditions. The half-logarithmic result of~\cite{LGX} controls a window whose lower endpoint has expected solution count $2^{\Theta(\sqrt N)}$. Our criterion concerns a bounded expected count. Concentration in the former window can therefore coexist with failure of Poisson counting in the latter.

Within the constant-mean window, the parameter $a_N$ measures the contribution of local neighbors. For $a_N\to0$, the central Hamming shells determine the second factorial moment. At $a_N\to a>0$, local neighbors contribute $\lambda(e^a-1)$. If $a_N\to\infty$, a conditioned solution has a diverging number of neighbors, whereas an unconditioned formula is unsatisfiable with high probability. Rare satisfiable formulas then carry the positive limiting mean.

The redundancy theorem explains why this change in counting law need not prevent solution blocking. Its first-moment argument requires a positive lower bound on the uniqueness probability, which the coexistence theorem provides even at bounded width offsets. The resulting formulas have $\Theta(N^2)$ clauses, and all but $o(M)$ slots are redundant with high probability conditional on uniqueness.

\subsection{The boundary distribution}
The main open problem is to determine the law of $X$ when $a_N\to a\in(0,\infty)$. Any proposed limit must account for the second factorial moment and the Poisson number of free directions seen from a conditioned solution. The local cube theorem describes a subset of the solution component. Identifying maximal components and controlling their joint occurrence would be needed to derive a law from component counts.

The same problem may also clarify coexistence for $\lambda\ge1$. Our uniqueness bound is obtained by subtracting the expected number of other solutions from the probability that the conditioned root has no free directions. It gives the factor $1-\lambda$ and hence a positive bound only for $0<\lambda<1$. Extending coexistence beyond this range requires finer control of other solutions around an isolated root. Quantitative error bounds would also help assess uniqueness probabilities at finite sizes.

\section{Conclusion}
The scale $N2^{-K}$ separates three forms of critical solution counting. Above logarithmic width, its convergence to zero is equivalent to both second-factorial-moment independence and an ordinary Poisson limit. At a positive finite limit, compatible free directions create local cubes and an explicit excess in the second factorial moment. When the scale diverges, the satisfiability probability tends to zero while the expected count stays positive.

For $0<\lambda<1$, unique satisfiability and unsatisfiability coexist exactly when $N2^{-K}$ remains bounded. Poisson convergence is therefore a stronger requirement than the one needed for coexistence and redundant-clause solution blocking. Determining the full boundary law would extend this description from moment and probability bounds to exact counting probabilities.

\end{document}